\documentclass{amsart}
\usepackage{amssymb}
\usepackage{amsmath}
\usepackage{amsfonts}

\newtheorem{theorem}{Theorem}[section]
\theoremstyle{plain}

\newtheorem{corollary}[theorem]{Corollary}
\newtheorem{definition}[theorem]{Definition}
\newtheorem{example}[theorem]{Example}
\newtheorem{lemma}[theorem]{Lemma}
\newtheorem{proposition}[theorem]{Proposition}
\newtheorem{remark}[theorem]{Remark}

\numberwithin{theorem} {section}
\input{tcilatex}

\begin{document}
\title[On frames in Hilbert modules over pro-$C^{\ast }$-algebras]{On frames
in Hilbert modules over pro-$C^{\ast }$-algebras}
\author{Maria Joi\c{t}a }
\address{Deparment of Mathematics, University of Bucharest, Bd. Regina
Elisabeta nr. 4-12, Bucharest, Romania }
\email{mjoita@fmi.unibuc.ro}
\urladdr{htp://www.geocities.com/m\_joita }
\thanks{This research was supported by CNCSIS (Romanian National Council for
Research in High Education) grant-code A 1065/2006.}

\begin{abstract}
We introduce the concept of frame of multipliers in Hilbert modules over pro-%
$C^{\ast }$-algebras and show that many properties of frames in Hilbert $%
C^{\ast }$-modules are valid for frames of multipliers in Hilbert modules
over pro-$C^{\ast }$-algebras.
\end{abstract}

\maketitle

\section{Introduction}

Hilbert $C^{\ast }$-modules are generalizations of Hilbert spaces by
allowing the inner product to take values in a $C^{\ast }$-algebra rather
than in the field of complex numbers. But the theory of Hilbert $C^{\ast }$%
-modules is different from the theory of Hilbert spaces, for example, no any
closed submodule of a Hilbert $C^{\ast }$-module is complemented. The notion
of frames in Hilbert $C^{\ast }$-modules was introduced and some properties
were investigated in \cite{2,3,13}. Using the well known Kasparov
stabilisation theorem \cite{10} which states that any countably generated
Hilbert $C^{\ast }$-module over $A$ is unitarily equivalent with a
complemented submodule of $H_{A}$, Frank and Larson \cite{2} showed that any
countably generated Hilbert module has a standard normalised frame. In \cite%
{13}, Raeburn and Thompson showed that the Kasparov stabilisation theorem is
valid for countably generated Hilbert $C^{\ast }$-modules in the multiplier.
Also they defined the concept of standard frame of multipliers for a Hilbert 
$C^{\ast }$-module, and showed that every Hilbert $C^{\ast }$-module
countably generated in the multiplier module admits a frame of multipliers,
thus generalizing results of Frank and Larson \cite{2}.

Pro-$C^{\ast }$-algebras are generalizations of $C^{\ast }$-algebras.
Instead of being given by a single $C^{\ast }$-norm, the topology of a pro-$%
C^{\ast }$-algebra is given by a directed family of $C^{\ast }$-seminorms.
In the literature, pro-$C^{\ast }$-algebras have been given by different
names such as $b^{\ast }$-algebras ( C. Apostol), $LMC^{\ast }$-algebras (G.
Lessner, K. Schm\"{u}dgen) or locally $C^{\ast }$-algebras (A. Inoue, M.\
Fragoulopoulou, A. Mallios, etc). Hilbert modules over pro-$C^{\ast }$%
-algebras were considered independently by Mallios and Phillips \cite{12}.
Phillips showed that the Kasparov stabilisation theorem is valid for
countably generated Hilbert modules over metrizable pro-$C^{\ast }$-algebras
and we shoved that this theorem is valid for countably generated Hilbert
modules over arbitrary pro-$C^{\ast }$-algebras \cite{7}.

In this paper we extend some results from \cite{2,3,13} in the context of
Hilbert modules over pro-$C^{\ast }$-algebras. The paper is organized as
follows. In Section 2 we recall some facts about pro-$C^{\ast }$-algebras
and Hilbert modules over pro-$C^{\ast }$-algebras. In Section 3, we
introduce the concept of frame of multipliers for Hilbert modules over pro-$%
C^{\ast }$-algebras and prove that any Hilbert module over a pro-$C^{\ast }$%
-algebra countably generated in the multiplier module admits a standard
normalised frame of multipliers. Also we show that the reconstruction
formula is valid for standard normalised frames of multipliers and the
existence of the frame transform. It is known that if the bounded part $b(E)$
of a Hilbert module $E$ over a pro-$C^{\ast }$-algebra $A$ is a Hilbert $%
C^{\ast }$-module over $b(A)$, and if $b(E)$ is countably generated in the
multiplier module, then $E$ is countably generated in the multiplier module.
We show that $b(E)$ admits a standard frame of multipliers, then $E$ admits
a standard frame of multipliers. In Section 4, we introduce the notion of
dual frame of multipliers and prove a necessary and sufficient condition for
that two frames of multipliers are duals to each other.

\section{Preliminaries}

Let $A$ be a pro-$C^{\ast }$-algebra. The set $S(A)$ of all continuous $%
C^{\ast }$-seminorms on $A$ is directed ($p\geq q$ if $p(a)\geq q(a)$ for
all $a\in A$). For each $p\in S(A)$, the quotient $\ast $-algebra $A/\ker p$%
, where $\ker p=\{a\in A;p(a)=0\}$, denoted by $A_{p}$, is a $C^{\ast }$%
-algebra in the $C^{\ast }$-norm induced by $p$. The canonical map from $A$
onto $A_{p}$ is denoted by $\pi _{p}^{A}$. For $p,q\in S(A)$ with $p\geq q$,
there is a canonical surjective morphism of $C^{\ast }$-algebras $\pi
_{pq}^{A}:A_{p}\rightarrow A_{q}$ such that $\pi _{pq}^{A}(\pi
_{p}^{A}(a))=\pi _{q}^{A}(a)$ for all $a\in A$, and $\{A_{p};\pi
_{pq}^{A}\}_{p,q\in S(A),p\geq q}$ is an inverse system of $C^{\ast }$%
-algebras. Moreover, the pro-$C^{\ast }$-algebras $A$ and $\lim\limits_{%
\underset{p\in S(A)}{\leftarrow }}A_{p}$ can be identified.

An element $a\in A$ is bounded if 
\begin{equation*}
\left\Vert a\right\Vert _{\infty }=\sup \{p(a);p\in S(A)\}<\infty .
\end{equation*}%
The set $b(A)$ of all bounded elements in $A$ is a $C^{\ast }$-algebra in
the $C^{\ast }$-norm $\left\Vert \cdot \right\Vert _{\infty }$. Moreover, $%
b(A)$ is dense in $A$.

Here we recall some facts about Hilbert modules over pro-$C^{\ast }$%
-algebras from \cite{7,12}.

\begin{definition}
A pre-Hilbert$\ A$-module is a complex vector space$\ E$\ which is also a
right $A$-module, compatible with the complex algebra structure, equipped
with an $A$-valued inner product $\left\langle \cdot ,\cdot \right\rangle
_{E}:E\times E\rightarrow A\;$which is $\mathbb{C}$-and $A$-linear in its
second variable and satisfies the following relations:

\begin{enumerate}
\item $\left\langle \xi ,\eta \right\rangle _{E}^{\ast }=\left\langle \eta
,\xi \right\rangle _{E}\;\;$for every $\xi ,\eta \in E;$

\item $\left\langle \xi ,\xi \right\rangle _{E}\geq 0\;\;$for every $\xi \in
E;$

\item $\left\langle \xi ,\xi \right\rangle _{E}=0\;$\ if and only if $\xi
=0. $
\end{enumerate}

We say that $E\;$is a Hilbert $A$-module if $E\;$is complete with respect to
the topology determined by the family of seminorms $\{\overline{p}%
_{E}\}_{p\in S(A)}\;$where $\overline{p}_{E}(\xi )=\sqrt{p\left(
\left\langle \xi ,\xi \right\rangle _{E}\right) },\xi \in E$.\smallskip
\end{definition}

An element $\xi $ in a Hilbert $A$-module $E$ is bounded if 
\begin{equation*}
\left\Vert \xi \right\Vert _{\infty }=\sup \{\overline{p}_{E}(\xi );p\in
S(A)\}<\infty .
\end{equation*}%
The set $b(E)$ of all bounded elements is a Hilbert $b(A)$-module which is
dense in $E$.

A Hilbert $A$-module $E$ is countably generated if there is a countable set $%
\{\xi _{n}\}_{n}$ in $E$ such that the submodule of $E$ generated by $\{\xi
_{n}a;a\in A,n=1,2,...\}$ is the whole of $E$.

If $A$ is a pro-$C^{\ast }$-algebra, then $A$ is a Hilbert $A$-module with $%
\left\langle a,b\right\rangle _{A}=a^{\ast }b,$ and the set $H_{A}$ of all
sequences $(a_{n})_{n}$ with $a_{n}\in A$ such that $\tsum\limits_{n}a_{n}^{%
\ast }a_{n}$ converges in $A$ is a Hilbert $A$-module with the action of $A$
on $H_{A}$ defined by $(a_{n})_{n}b=(a_{n}b)_{n}$ and the inner product
defined by $\left\langle (a_{n})_{n},(b_{n})_{n}\right\rangle
_{H_{A}}=\tsum\limits_{n}a_{n}^{\ast }b_{n}.$ Moreover, if $A$ has a
countable approximate unit, then the Hilbert $A$-modules $A$ and $H_{A}$ are
countably generated.

Let $E\;$be a Hilbert $A$-module.\ For $p\in S(A)$,$\;$the quotient vector
space $E/\ker \left( \overline{p}_{E}\right) ,$ where $\ker \left( \overline{%
p}_{E}\right) =\{\xi \in E;\overline{p}_{E}(\xi )=0\}$,$\;$denoted by $E_{p}$%
, is a Hilbert $A_{p}$-module with $(\xi +\ker \left( \overline{p}%
_{E}\right) {})\pi _{p}^{A}(a)=\xi a+\ker \left( \overline{p}_{E}\right)
{}\; $and $\left\langle \xi +\ker \left( \overline{p}_{E}\right) {},\eta
+\ker \left( \overline{p}_{E}\right) {}\right\rangle _{E_{p}}=\pi
_{p}^{A}(\left\langle \xi ,\eta \right\rangle _{E})$.\ The canonical map
from $E\;$onto $E_{p}$ is denoted by $\sigma _{p}^{E}$. For $p,q\in S(A)$
with $p\geq q\;$there is a canonical surjective morphism of vector spaces $%
\sigma _{pq}^{E}\;$from $E_{p}\;$onto $E_{q}\;$such that $\sigma
_{pq}^{E}(\sigma _{p}^{E}(\xi ))=\sigma _{q}^{E}(\xi )$ for all $\xi \in
E,\; $and $\ \{E_{p};A_{p};\sigma _{pq}^{E},\pi _{pq}^{A}$\ $\}_{p,q\in
S(A),p\geq q}$ is an inverse system of Hilbert $C^{\ast }$-modules in the
following sense: $\sigma _{pq}^{E}(\xi _{p}a_{p})=\sigma _{pq}^{E}(\xi
_{p})\pi _{pq}^{A}(a_{p}),\xi _{p}\in E_{p},a_{p}\in A_{p};$ $\left\langle
\sigma _{pq}^{E}(\xi _{p}),\sigma _{pq}^{E}(\eta _{p})\right\rangle
_{E_{q}}=\pi _{pq}^{A}(\left\langle \xi _{p},\eta _{p}\right\rangle
_{E_{p}}),$\ $\ $\ $\xi _{p},\eta _{p}\in E_{p};$ $\sigma _{pp}^{E}(\xi
_{p})=\xi _{p},\;\xi _{p}\in E_{p}\;$and $\sigma _{qr}^{E}\circ \sigma
_{pq}^{E}=\sigma _{pr}^{E}\;$if $p\geq q\geq r$.$\ $The Hilbert $A$-modules $%
\lim\limits_{\underset{p\in S(A)}{\leftarrow }}E_{p}$ and $E$ can be
identified.

Given two Hilbert $A$-modules $E$ and $F$, a module morphism $T:E\rightarrow
F\;$is continuous if for each $p\in S(A)$ there is $M_{p}>0$ such that $%
\overline{p}_{F}\left( T\xi \right) \leq M_{p}\overline{p}_{E}\left( \xi
\right) $ for all $\xi \in E$, and it is adjointable if there is a module
morphism $T^{\ast }:F\rightarrow E$\ such that $\left\langle T\xi ,\eta
\right\rangle =\left\langle \xi ,T^{\ast }\eta \right\rangle \;$for every $%
\xi \in E$ and $\eta \in F$. Any adjointable module morphism is continuous.
The set of all adjointable module morphisms from $E$ to $F$ is denoted by $%
L\left( E,F\right) $. For each $p\in S(A),$ there is a linear map $\left(
\pi _{p}^{E,F}\right) _{\ast }:L(E,F)\rightarrow L\left( E_{p},F_{p}\right) $
defined by 
\begin{equation*}
\left( \pi _{p}^{E,F}\right) _{\ast }\left( T\right) \left( \sigma
_{p}^{E}(\xi )\right) =\sigma _{p}^{F}\left( T\left( \xi \right) \right)
\end{equation*}%
for $T\in L(E)$ and $\xi \in E$. The vector space $L(E,F)$ is a complete
locally convex space with respect to the topology defined by the family of
seminorms $\{\widetilde{p}_{L(E,F)}\}_{p\in S(A)}$, where $\widetilde{p}%
_{L(E,F)}(T)=\left\Vert (\pi _{p}^{E,F})_{\ast }(T)\right\Vert
_{L(E_{p},F_{p})},$ $T\in L(E,F)$. If $E=F,$ $L(E,E)$ is denoted by $L(E)$
and it is a pro-$C^{\ast }$-algebra.

For $p,q\in S(A)$ with $p\geq q$, there is a linear map $\left( \pi
_{pq}^{E,F}\right) _{\ast }:L(E_{p},F_{p})\rightarrow L(E_{q},F_{q})$ such
that 
\begin{equation*}
\left( \pi _{pq}^{E,F}\right) _{\ast }\left( T_{p}\right) \left( \sigma
_{q}^{E}\left( \xi \right) \right) =\sigma _{pq}^{F}\left( T_{p}\left(
\sigma _{p}^{E}\left( \xi \right) \right) \right)
\end{equation*}
for all $T_{p}\in L(E_{p},F_{p})$ and $\xi \in E$, and $\{L(E_{p},F_{p});$ $%
\;(\pi _{pq}^{E,F})_{\ast }\}_{p,q\in S(A),p\geq q}$ is an inverse system of
Banach spaces. Moreover, the complete locally convex spaces $\lim\limits_{%
\underset{p\in S(A)}{\leftarrow }}L(E_{p},F_{p})$ and $L(E,F)$ are
isomorphic.

The set $b(L(E,F))$ of all bounded elements in $L(E,F)$ ($T\in b(L(E,F)$ if $%
\sup \{$\ $\ \widetilde{p}_{L(E,F)}(T);p\in S(A)\}<\infty $), is a Banach
space which is isomorphic with $L(b(E),b(F))$.

Now we recall some facts about multiplier modules from \cite{8,13}. The set $%
L(A,E)$ of all adjointable module morphisms from $A$ to $E$ is a Hilbert $%
L(A)$-module with the action of $L(A)$ on $L(A,E)$ defined by 
\begin{equation*}
L(A,E)\times L(A)\backepsilon \left( T,S\right) \mapsto T\cdot S=T\circ S\in
L(A,E)
\end{equation*}%
and the inner-product defined by 
\begin{equation*}
L(A,E)\times L(A,E)\backepsilon \left( T,R\right) \mapsto \left\langle
T,S\right\rangle =T^{\ast }\circ S\in L(A).
\end{equation*}%
Since the locally $C^{\ast }$-algebras $L(A)$ and $M(A),$ the multiplier
algebra of $A$, can be identified \cite{12,7}, the Hilbert $L(A)$-module $%
L(A,E)$ can be regarded as a Hilbert $M(A)$-module. The Hilbert $M(A)$%
-module $L(A,E)$ is called \textit{the multiplier module} of $E$ and it is
denoted by $M(E)$. Moreover, the topology on $M(E)$ induced by the inner
product coincides with the topology defined by the family of seminorms $\{%
\overline{p}_{M(E)}\}_{p\in S(A)}$, with 
\begin{equation*}
\overline{p}_{M(E)}(h)=\widetilde{p}_{L(A,E)}\left( h\right)
\end{equation*}%
for all $h\in M(E)\;$and for all $p\in S(A)$. The map i$_{E}:E\rightarrow
M(E)$ defined by 
\begin{equation*}
\text{i}_{E}(\xi )\left( a\right) =\xi a,\xi \in E,a\in A
\end{equation*}%
identifies $E$ with a Hilbert submodule of $M(E)$ and then 
\begin{equation*}
\left\langle h,\xi \right\rangle _{M(E)}=h^{\ast }\left( \xi \right)
\end{equation*}%
for all $h\in M(E)$ and $\xi \in E$. Moreover, if $a\in A$ and $h\in M(E)$,
then $h\cdot a$ can be identified with $h(a)$.

\begin{definition}
A Hilbert $A$-module $E$ is countably generated in $M(E)$ if there is a
countable set $\{h_{n};$ $h_{n}\in M(E),$ $n=1,2,...\}$ such that the closed
submodule of $M(E)$ generated by $\{h_{n}\cdot a;$ $a\in A,$ $n=1,2,...\}$
is the whole of $E$.
\end{definition}

\begin{remark}
If the Hilbert $A$-module $E$ is countably generated in $M(E)$, then, for
each $p\in S(A),$ $\{(\pi _{p}^{A,E})_{\ast }(h_{n});$ $h_{n}\in M(E),$ $%
n=1,2,...\}\subseteq M(E_{p})$ is a generating set for $E_{p}$, since $(\pi
_{p}^{A,E})_{\ast }(h_{n})\cdot \pi _{p}^{A}(a)=\sigma _{p}^{E}\left(
h_{n}\cdot a\right) $ for all $n=1,2,...$ and for all $a\in A$, and since $%
\sigma _{p}^{E}(E)=E_{p}$. Therefore the Hilbert $A_{p}$-module $E_{p}$ is
countably generated in $M(E_{p})$ for each $p\in S(A)$.
\end{remark}

\begin{example}
If $A$ is a pro-$C^{\ast }$-algebra, then the Hilbert $A$-module $A$ is
countably generated in $M(A)$, $\{1_{M(A)}\}$ being a generating set.
\end{example}

\begin{example}
For any pro-$C^{\ast }$-algebra $A$, the Hilbert $A$-module $H_{A}$ is
countably generated in $M(H_{A})$.

Indeed, for each positive integer $n$ consider the linear map $%
e_{n}:A\rightarrow H_{A}$ by $e_{n}(a)=(0,...0,a,0,...)$, the element in $%
H_{A}$ whose all the components are $0$ excepts at $n^{\text{th}}$ component
which is $a$. Clearly, $e_{n}$ is a module morphism. Moreover, $e_{n}$ is
adjointable and $e_{n}^{\ast }\left( \left( a_{m}\right) _{m}\right) =$ $%
a_{n}$.

Let $\left( a_{n}\right) _{n}\in H_{A}$. Since 
\begin{equation*}
\overline{p}_{H_{A}}\left( \left( a_{n}\right)
_{n}-\sum\limits_{k=1}^{m}e_{k}\cdot a_{k}\right) ^{2}=p\left(
\sum\limits_{k=m+1}^{\infty }a_{k}^{\ast }a_{k}\right)
\end{equation*}%
for all $p\in S(A)$ and for all positive integer $m,$ $\left( a_{n}\right)
_{n}=\sum\limits_{k}e_{k}\cdot a_{k}.$ Therefore, the Hilbert submodule of $%
M(H_{A})$ generated by $\{e_{m}\cdot a;a\in A,m=1,2,...\}$ is $H_{A}$, and
so $\{e_{m},m=1,2,...\}$ $\subseteq M(H_{A})$ is a generating set for $H_{A}$%
.
\end{example}

\begin{remark}
If $E$ is a countably generated Hilbert $A$-module, then $E$ is countably
generated in $M(E)$.

In general, $E$ is not countably generated when $E$ is countably generated
in $M(E)$.

\textbf{Example.} Let $A$ be a pro-$C^{\ast }$-algebra which does not have a
countable approximate unit. We seen that the Hilbert $A$-module $A$ is
countably generated in $M(A)$ but it is not countably generated.
\end{remark}

\begin{remark}
If $E$ is a Hilbert $A$-module such that $b(E)$ is countably generated in $%
M(b(E)),$ then $E$ is countably generated in $M(E)$.
\end{remark}

\section{Frame transform and reconstruction frame}

Let $E$ be a Hilbert $A$-module.

\begin{definition}
A sequence $\{h_{n}\}_{n}$ in $M(E)$ is a standard frame of multipliers in $%
E $ if for each $\xi \in E$, $\sum\limits_{n}\left\langle \xi
,h_{n}\right\rangle _{M(E)}\left\langle h_{n},\xi \right\rangle _{M(E)}$
converges in $A$, and there are two positive constants $C$ and $D$ such that 
\begin{equation*}
C\left\langle \xi ,\xi \right\rangle _{E}\leq \sum\limits_{n}\left\langle
\xi ,h_{n}\right\rangle _{M(E)}\left\langle h_{n},\xi \right\rangle
_{M(E)}\leq D\left\langle \xi ,\xi \right\rangle _{E}
\end{equation*}%
for all $\xi \in E$. If $D=C=1$ we say that $\{h_{n}\}_{n}$ is a standard
normalized frame of multipliers.
\end{definition}

\begin{remark}
Let $\{h_{n}\}_{n}$ be a sequence in $M(E)$. If $\{h_{n}\}_{n}$ is a
standard frame of multipliers in $E$, then $\{\left( \pi _{p}^{A,E}\right)
_{\ast }\left( h_{n}\right) \}_{n}$ is a standard frame of multipliers in $%
E_{p}$ for each $p\in S(A)$.
\end{remark}

\begin{remark}
Let $\{h_{n}\}_{n}$ be a sequence in $M(E)$. Then $\{h_{n}\}_{n}$ is a
standard normalised frame of multipliers in $E$ if and only if $\{\left( \pi
_{p}^{A,E}\right) _{\ast }\left( h_{n}\right) \}_{n}$ is a standard
normalised frame of multipliers in $E_{p}$ for each $p\in S(A).$
\end{remark}

\begin{remark}
If $\{h_{n}\}_{n}$ is a standard frame of multipliers in $E$, then $h_{n}\in
b(M(E))$ for all positive integer $n$.

Indeed, let $\xi \in b(E)$ and let $m$ be a positive integer. From 
\begin{equation*}
0\leq \left\langle \xi ,h_{m}\right\rangle _{M(E)}\left\langle h_{m},\xi
\right\rangle _{M(E)}\leq \sum\limits_{n}\left\langle \xi
,h_{n}\right\rangle _{M(E)}\left\langle h_{n},\xi \right\rangle _{M(E)}\leq
D\left\langle \xi ,\xi \right\rangle _{E}
\end{equation*}%
and \cite{4}, we deduce that 
\begin{equation*}
p\left( \left\langle h_{m},\xi \right\rangle _{M(E)}\right) ^{2}\leq
Dp\left( \left\langle \xi ,\xi \right\rangle _{E}\right)
\end{equation*}%
for all $p\in S(A)$. Then 
\begin{equation*}
\overline{p}_{A}\left( h_{m}^{\ast }\left( \xi \right) \right) ^{2}\leq D%
\overline{p}_{E}\left( \xi \right) ^{2}\leq D\left\Vert \xi \right\Vert
_{\infty }^{2}
\end{equation*}%
for all $p\in S(A)$. This implies that $h_{m}^{\ast }\in b(M(E))$. Therefore 
$h_{m}\in b(M(E))$.
\end{remark}

\begin{example}
For any pro-$C^{\ast }$-algebra, $\{e_{n}\}_{n}$ is a standard normalised
frame of multipliers in $H_{A}$. Indeed, if $\left( a_{n}\right) _{n}\in
H_{A}$ then, since 
\begin{equation*}
\left\langle \left( a_{m}\right) _{m},e_{n}\right\rangle _{M\left(
H_{A}\right) }\left\langle e_{n},\left( a_{m}\right) _{m}\right\rangle
_{M\left( H_{A}\right) }=a_{n}^{\ast }a_{n}
\end{equation*}%
for each positive integer $n$, we have 
\begin{equation*}
\sum\limits_{n}\left\langle \left( a_{m}\right) _{m},e_{n}\right\rangle
_{M\left( H_{A}\right) }\left\langle e_{n},\left( a_{m}\right)
_{m}\right\rangle _{M\left( H_{A}\right) }=\sum\limits_{n}a_{n}^{\ast
}a_{n}=\left\langle \left( a_{m}\right) _{m},\left( a_{m}\right)
_{m}\right\rangle _{H_{A}}
\end{equation*}%
and so $\{e_{n}\}_{n}$ is a standard normalised frame of multipliers in $%
H_{A}$.
\end{example}

\begin{proposition}
Any countably generated Hilbert $A$-module $E$ in $M(E)$ admits a standard
normalised frame of multipliers.
\end{proposition}

\proof%
Indeed, let $P:H_{A}\rightarrow E$ be the projection of $H_{A}$ on $E$ \cite[%
Theorem 4.2]{8}. Then $\{\left( \pi _{p}^{H_{A},E}\right) _{\ast }(P)\circ
\left( \pi _{p}^{A,H_{A}}\right) _{\ast }\left( e_{n}\right) \}_{n}$ is a
standard normalised frame of multipliers in $E_{p}$ for each $p\in S(A)$ 
\cite[Corollary 3.3]{13}. From this fact, Remark 3.3 and taking into account
that $\left( \pi _{p}^{H_{A},E}\right) _{\ast }(P)\circ \left( \pi
_{p}^{A,H_{A}}\right) _{\ast }\left( e_{n}\right) =$ $\left( \pi
_{p}^{A,E}\right) _{\ast }\left( P\circ e_{n}\right) $ for all $n$ and for
all $p\in S(A)$ we conclude that $\{P\circ e_{n}\}_{n}$ is a standard
normalised frame of multipliers in $E$.%
\endproof%

\begin{theorem}
( The reconstruction formula) Let $E$ be a countably generated Hilbert $A$%
-module in $M(E)$ and let $\{h_{n}\}_{n}$ be a sequence in $M(E)$. Then $%
\{h_{n}\}_{n}$ is a standard normalised frame of multipliers if and only if
for all $\xi \in E$, $\sum\limits_{n}h_{n}\cdot \left\langle h_{n},\xi
\right\rangle _{M(E)}$ converges in $E$ and moreover, 
\begin{equation*}
\xi =\sum\limits_{n}h_{n}\cdot \left\langle h_{n},\xi \right\rangle _{M(E)}%
\text{.}
\end{equation*}
\end{theorem}

\proof
By Remark 3.3 and \cite[Theorem 3.4]{13}, $\{h_{n}\}_{n}$ is a standard
normalised frame of multipliers in $E$ if and only if $\sum\limits_{n}\left(
\pi _{p}^{A,E}\right) _{\ast }\left( h_{n}\right) \cdot \left\langle \left(
\pi _{p}^{A,E}\right) _{\ast }\left( h_{n}\right) ,\sigma _{p}^{E}\left( \xi
\right) \right\rangle _{M(E)}$ converges in $E_{p}$ for all $\xi \in E$ and
for each $p\in S(A)$, and moreover, 
\begin{equation*}
\sigma _{p}^{E}\left( \xi \right) =\sum\limits_{n}\left( \pi
_{p}^{A,E}\right) _{\ast }\left( h_{n}\right) \cdot \left\langle \left( \pi
_{p}^{A,E}\right) _{\ast }\left( h_{n}\right) ,\sigma _{p}^{E}\left( \xi
\right) \right\rangle _{M(E)}
\end{equation*}%
From this fact and taking into account that 
\begin{eqnarray*}
&&\overline{p}_{E}\left( \xi -\sum\limits_{k=1}^{n}h_{k}\cdot \left\langle
h_{k},\xi \right\rangle _{M(E)}\right) \\
&=&\left\Vert \sigma _{p}^{E}\left( \xi \right) -\sigma _{p}^{E}\left(
\sum\limits_{k=1}^{n}h_{k}\cdot \left\langle h_{k},\xi \right\rangle
_{M(E)}\right) \right\Vert _{E_{p}} \\
&=&\left\Vert \sigma _{p}^{E}\left( \xi \right) -\sum\limits_{k=1}^{n}\left(
\pi _{p}^{A,E}\right) _{\ast }\left( h_{k}\right) \cdot \left\langle \left(
\pi _{p}^{A,E}\right) _{\ast }\left( h_{k}\right) ,\sigma _{p}^{E}\left( \xi
\right) \right\rangle _{M(E)}\right\Vert _{E_{p}}
\end{eqnarray*}%
for all $\xi \in E$, for all $p\in S(A)$ and for all positive integer $n$,
we deduce that $\{h_{n}\}_{n}$ is a standard normalised frame of multipliers
in $E$ if and only if $\sum\limits_{n}h_{n}\cdot \left\langle h_{n},\xi
\right\rangle _{M(E)}$ converges in $E$ for all $\xi \in E$, and moreover, $%
\xi =\sum\limits_{n}h_{n}\cdot \left\langle h_{n},\xi \right\rangle _{M(E)}$
for all $\xi \in E$.%
\endproof%

\begin{remark}
If $\{h_{n}\}_{n}$ is a standard normalised frame of multipliers in $E$,
then\ $\tsum\limits_{n}(h_{n}\circ $ $h_{n}^{\ast })\left( \xi \right) =\xi $
for all $\xi \in E$, since $\left( h_{n}\circ h_{n}^{\ast }\right) \left(
\xi \right) =h_{n}\cdot \left\langle h_{n},\xi \right\rangle _{M(E)}$ for
each positive integer $n$. Therefore, $\{h_{n}\}_{n}$ is a standard
normalied frame of multipliers in $E$ if and only if $\tsum\limits_{n}\left(
h_{n}\circ h_{n}^{\ast }\right) \left( \xi \right) $ converges in $E$ for
each $\xi \in E$ and moreover, $\tsum\limits_{n}\left( h_{n}\circ
h_{n}^{\ast }\right) \left( \xi \right) =\xi $.
\end{remark}

\begin{corollary}
If $b(E)$ admits a standard normalised frame of multipliers, then $E$ admits
a standard normalised frame of multipliers.
\end{corollary}

\proof
We will show that if $\{h_{n}\}_{n}$ is a standard normalised frame of
multipliers in $b(E)$, then $\{\widetilde{h_{n}}\}_{n}$, where $\widetilde{%
h_{n}}$ is the extension of $h_{n}$ to an element in $M(E)$, is a standard
normalised frame of multipliers in $E$.

Let $\xi \in E$, $p\in S(A)$ and $\varepsilon >0$. Since $b(E)$ is dense in $%
E$, there is $\xi _{0}\in b(E)$ such that $\overline{p}_{E}\left( \xi -\xi
_{0}\right) \leq \varepsilon /3$. Since $\{h_{n}\}_{n}$ is a standard
normalised frame of multipliers in $b(E)$, there is $n_{0}$ such that 
\begin{equation*}
\left\Vert \xi _{0}-\tsum\limits_{k=1}^{n}\left( h_{k}\circ h_{k}^{\ast
}\right) \left( \xi _{0}\right) \right\Vert _{\infty }\leq \varepsilon /3
\end{equation*}%
for all $n$ with $n\geq n_{0}$. Then 
\begin{eqnarray*}
\overline{p}_{E}\left( \xi -\tsum\limits_{k=1}^{n}\left( \widetilde{h_{k}}%
\circ \widetilde{h_{k}}^{\ast }\right) \left( \xi \right) \right) &\leq &%
\overline{p}_{E}\left( \xi -\xi _{0}\right) +\overline{p}_{E}\left( \xi
_{0}-\tsum\limits_{k=1}^{n}\left( h_{k}\circ h_{k}^{\ast }\right) \left( \xi
_{0}\right) \right) \\
&&+\overline{p}_{E}\left( \tsum\limits_{k=1}^{n}\left( \widetilde{h_{k}}%
\circ \widetilde{h_{k}}^{\ast }\right) \left( \xi -\xi _{0}\right) \right) \\
&\leq &\varepsilon /3+\left\Vert \xi _{0}-\tsum\limits_{k=1}^{n}\left(
h_{k}\circ h_{k}^{\ast }\right) \left( \xi _{0}\right) \right\Vert _{\infty }
\\
&&+\widetilde{p}_{L(E)}\left( \tsum\limits_{k=1}^{n}\widetilde{h_{k}}\circ 
\widetilde{h_{k}}^{\ast }\right) \overline{p}_{E}\left( \xi -\xi _{0}\right)
\\
&\leq &\varepsilon /3\left( 2+\left\Vert \tsum\limits_{k=1}^{n}\widetilde{%
h_{k}}\circ \widetilde{h_{k}}^{\ast }\right\Vert _{\infty }\right) \\
&=&\varepsilon /3\left( 2+\left\Vert \tsum\limits_{k=1}^{n}h_{k}\circ
h_{k}^{\ast }\right\Vert \right) \\
&\leq &\varepsilon
\end{eqnarray*}%
for all $n$ with $n\geq n_{0}.$ This shows that $\tsum\limits_{n}\left( 
\widetilde{h_{n}}\circ \widetilde{h_{n}}^{\ast }\right) \left( \xi \right) $
converges to $\xi $ in $E$ for each $\xi \in E$ and so $\{\widetilde{h_{n}}%
\}_{n}$ is a standard normalised frame of multipliers in $E.$ 
\endproof%

If $\{h_{n}\}_{n}$ is a standard frame of multipliers in $E$, then $%
\sum\limits_{n}\left\langle \xi ,h_{n}\right\rangle _{M(E)}\left\langle
h_{n},\xi \right\rangle _{M(E)}$ converges in $E$ for all $\xi \in E$. From
this fact and taking into account that $\left\langle h_{n},\xi \right\rangle
_{M(E)}\in A$ for all positive integer $n$, we conclude that $\left(
\left\langle h_{n},\xi \right\rangle _{M\left( E\right) }\right) _{n}\in
H_{A}$. Thus we can define a linear map $\theta :E\rightarrow H_{A}$ by 
\begin{equation*}
\theta \left( \xi \right) =\left( \left\langle h_{n},\xi \right\rangle
_{M\left( E\right) }\right) _{n}.
\end{equation*}%
Moreover, $\theta $ is a continuous module morphism, since 
\begin{equation*}
\theta \left( \xi a\right) =\left( \left\langle h_{n},\xi a\right\rangle
_{M\left( E\right) }\right) _{n}=\left( \left\langle h_{n},\xi \right\rangle
_{M\left( E\right) }a\right) _{n}=\theta \left( \xi \right) a
\end{equation*}%
for all $\xi \in E$ and for all $a\in A$ and 
\begin{equation*}
\overline{p}_{H_{A}}\left( \theta \left( \xi \right) \right) ^{2}=p\left(
\sum\limits_{n}\left\langle \xi ,h_{n}\right\rangle _{M(E)}\left\langle
h_{n},\xi \right\rangle _{M(E)}\right) \leq C\overline{p}_{E}\left( \xi
\right) ^{2}
\end{equation*}%
for all $\xi \in E$ and for all $p\in S(A)$.

\begin{definition}
Let $\{h_{n}\}_{n}$ be a standard frame of multipliers in $E$. The module
morphism $\theta :E\rightarrow H_{A}$ defined by $\theta \left( \xi \right)
=\left( \left\langle h_{n},\xi \right\rangle _{M\left( E\right) }\right)
_{n} $ is called the frame transform for $\{h_{n}\}_{n}$.
\end{definition}

\begin{theorem}
( The frame transform) Let $E$ be a countably generated Hilbert $A$-module
in $M(E)$ and let $\{h_{n}\}_{n}$ be a standard frame of multipliers in $E$.
The frame transform $\theta $ is an adjointable module morphism which
realizes an embedding of $E$ onto an orthogonal summand of $H_{A}$, and $%
\theta ^{\ast }\circ e_{n}=h_{n}$ for all $n$. Moreover, $\theta ^{\ast
}\circ \theta $ is an invertible element in $b(L\left( E\right) )$.
\end{theorem}

\proof%
Since $\{h_{n}\}_{n}$ is a standard frame of multipliers in $E$, $\{\left(
\left( \pi _{p}^{A,E}\right) _{\ast }\left( h_{n}\right) \right) \}_{n}$ is
a standard frame of multipliers in $E_{p}$ for each $p\in S(A)$. Let $p\in
S(A)$. The frame transform for $\{\left( \left( \pi _{p}^{A,E}\right) _{\ast
}\left( h_{n}\right) \right) \}_{n}$ is an adjointable operator $\theta
_{p}:E_{p}\rightarrow H_{A_{p}}$ defined by 
\begin{equation*}
\theta _{p}\left( \sigma _{p}^{E}\left( \xi \right) \right) =\left(
\left\langle \left( \pi _{p}^{A,E}\right) _{\ast }\left( h_{n}\right)
,\sigma _{p}^{E}\left( \xi \right) \right\rangle _{M(E_{p})}\right) _{n}.
\end{equation*}%
Moreover, $\theta _{p}$ preserves the inner product and $\theta _{p}^{\ast
}\circ \left( \pi _{p}^{A,H_{A}}\right) _{\ast }\left( e_{n}\right) =\left(
\pi _{p}^{A,E}\right) _{\ast }\left( h_{n}\right) $ \cite[Theorem 3.5]{13}.
Since 
\begin{equation*}
\left( \pi _{pq}^{E,H_{A}}\right) _{\ast }\left( \theta _{p}\right) =\theta
_{q}
\end{equation*}%
for all $p,q\in S(A)$ with $p\geq q$, there is $\theta \in L(E,H_{A})$ such
that 
\begin{equation*}
\left( \pi _{p}^{E,H_{A}}\right) _{\ast }\left( \theta \right) =\theta _{p}
\end{equation*}%
for all $p\in S(A)$. Clearly 
\begin{equation*}
\theta \left( \xi \right) =\left( \left\langle h_{n},\xi \right\rangle
_{M\left( E\right) }\right) _{n}
\end{equation*}%
for all $\xi \in E$. Therefore $\theta $ is the frame transform for $%
\{h_{n}\}_{n}$. Moreover, $\theta $ preserves the inner product and $\theta
^{\ast }\circ e_{n}=h_{n}$ for all $n$.

From 
\begin{eqnarray*}
C\overline{p}_{E}\left( \xi \right) ^{2} &=&Cp\left( \left\langle \xi ,\xi
\right\rangle _{E}\right) \leq p\left( \sum\limits_{n}\left\langle \xi
,h_{n}\right\rangle _{M(E)}\left\langle h_{n},\xi \right\rangle
_{M(E)}\right) \\
&=&\overline{p}_{H_{A}}\left( \theta \left( \xi \right) \right) ^{2}\leq
Dp\left( \left\langle \xi ,\xi \right\rangle _{E}\right) =D\overline{p}%
_{E}\left( \xi \right) ^{2}
\end{eqnarray*}%
for all $p\in S(A)$ and for all $\xi \in E$, we conclude that $\theta $ is
an injective adjointable module morphism from $E$ to $H_{A}$ with closed
range. Moreover, $\theta \in b(L(E,H_{A})$. By \cite[Theorem 2.2]{5}, $%
\theta \left( E\right) $ has an orthogonal complement in $H_{A},$ $\theta
^{\ast }$ is surjective and the restriction $\theta ^{\ast }|_{\theta (E)}$
of $\theta ^{\ast }$ on $\theta (E)$ is an invertible element in $b(L(\theta
(E),H_{A})$. Then $\theta ^{\ast }\circ \theta $ is invertible, and
moreover, $\theta ^{\ast }\circ \theta \in b(L(E))$. 
\endproof%

\begin{theorem}
Let $E$ be a countably generated Hilbert $A$-module in $M(E)$ and let $%
\{h_{n}\}_{n}$ be a sequence in $M(E)$. Then $\{h_{n}\}_{n}$ is a standard
frame of multipliers in $E$ if and only if there is an invertible element $T$
in $b(L(E))$ such that $\{T\circ h_{n}\}_{n}$ is a standard normalised frame
of multipliers in $E$.
\end{theorem}

\proof%
Suppose that $\{h_{n}\}_{n}$ is a standard frame of multipliers in $E$. Let $%
\theta $ be the frame transform. Then $\theta ^{\ast }\circ \theta $, and so 
$\left( \theta ^{\ast }\circ \theta \right) ^{-\frac{1}{2}}$, is an
invertible element in $b(L(E))$. Let $\xi \in E$. Then there is $\eta \in E$
such that $\left( \theta ^{\ast }\circ \theta \right) ^{\frac{1}{2}}\left(
\eta \right) =\xi $ and 
\begin{eqnarray*}
\left\langle \xi ,\xi \right\rangle _{E} &=&\left\langle \left( \theta
^{\ast }\circ \theta \right) ^{\frac{1}{2}}\left( \eta \right) ,\left(
\theta ^{\ast }\circ \theta \right) ^{\frac{1}{2}}\left( \eta \right)
\right\rangle _{E}=\left\langle \theta \left( \eta \right) ,\theta \left(
\eta \right) \right\rangle _{H_{A}} \\
&=&\sum\limits_{n}\left\langle \eta ,h_{n}\right\rangle _{M(E)}\left\langle
h_{n},\eta \right\rangle _{M(E)} \\
&=&\sum\limits_{n}\left\langle \left( \theta ^{\ast }\circ \theta \right) ^{-%
\frac{1}{2}}(\xi ),h_{n}\right\rangle _{M(E)}\left\langle h_{n},\left(
\theta ^{\ast }\circ \theta \right) ^{-\frac{1}{2}}(\xi )\right\rangle
_{M(E)} \\
&=&\sum\limits_{n}\left\langle \xi ,\left( \theta ^{\ast }\circ \theta
\right) ^{-\frac{1}{2}}\circ h_{n}\right\rangle _{M(E)}\left\langle \left(
\theta ^{\ast }\circ \theta \right) ^{-\frac{1}{2}}\circ h_{n},\xi
\right\rangle _{M(E)}\text{.}
\end{eqnarray*}%
Therefore $\{\left( \theta ^{\ast }\circ \theta \right) ^{-\frac{1}{2}}\circ
h_{n}\}_{n}$ is a standard normalised frame of multipliers in $E$.

Conversely, suppose that there is an invertible element $T$ in $b(L(E))$
such that $\{T\circ h_{n}\}_{n}$ is a standard normalised frame of
multipliers in $E$. Since $T\in b(L(E))$,%
\begin{eqnarray*}
\pi _{p}^{A}\left( \left\langle T(\xi ),T(\xi )\right\rangle _{E}\right)
&=&\left\langle \left( \pi _{p}^{E,E}\right) _{\ast }(T)(\sigma
_{p}^{E}\left( \xi \right) ),\left( \pi _{p}^{E,E}\right) _{\ast }(T)(\sigma
_{p}^{E}\left( \xi \right) )\right\rangle _{E_{p}} \\
&&\text{( see, for example, \cite[Proposition 1.2]{11})} \\
&\leq &\left\Vert \left( \pi _{p}^{E,E}\right) _{\ast }(T)\right\Vert
^{2}\left\langle \sigma _{p}^{E}\left( \xi \right) ,\sigma _{p}^{E}\left(
\xi \right) \right\rangle _{E_{p}} \\
&=&\widetilde{p}_{L(E)}(T)^{2}\pi _{p}^{A}\left( \left\langle \xi ,\xi
\right\rangle _{E}\right) \\
&\leq &\left\Vert T\right\Vert _{\infty }^{2}\pi _{p}^{A}\left( \left\langle
\xi ,\xi \right\rangle _{E}\right)
\end{eqnarray*}%
for all $p\in S(A)$ and for all $\xi \in E,$ and then by \cite{4} 
\begin{equation*}
\left\langle T(\xi ),T(\xi )\right\rangle _{E}\leq \left\Vert T\right\Vert
_{\infty }^{2}\left\langle \xi ,\xi \right\rangle _{E}
\end{equation*}%
for all $\xi \in E.$ Then, since $T$ is invertible we have 
\begin{equation*}
\left\Vert T\right\Vert _{\infty }^{-2}\left\langle \xi ,\xi \right\rangle
_{E}\leq \left\langle \left( T^{\ast }\right) ^{-1}\left( \xi \right)
,\left( T^{\ast }\right) ^{-1}\left( \xi \right) \right\rangle _{E}\leq
\left\Vert T^{-1}\right\Vert _{\infty }^{2}\left\langle \xi ,\xi
\right\rangle _{E}
\end{equation*}%
for all $\xi \in E$. From these relations and taking into account that

\begin{eqnarray*}
&&\left\langle \left( T^{\ast }\right) ^{-1}\left( \xi \right) ,\left(
T^{\ast }\right) ^{-1}\left( \xi \right) \right\rangle _{E} \\
&&\{T\circ h_{n}\}_{n}\ \text{is a standard normalised frame of }E\text{ in }%
M(E) \\
&=&\sum\limits_{n}\left\langle \left( T^{\ast }\right) ^{-1}\left( \xi
\right) ,T\circ h_{n}\right\rangle _{M(E)}\left\langle T\circ h_{n},\left(
T^{\ast }\right) ^{-1}\left( \xi \right) \right\rangle _{M(E)} \\
&=&\lim_{n}\sum\limits_{k=1}^{n}\left\langle \left( T^{\ast }\right)
^{-1}\left( \xi \right) ,T\circ h_{k}\right\rangle _{M(E)}\left\langle
T\circ h_{k},\left( T^{\ast }\right) ^{-1}\left( \xi \right) \right\rangle
_{M(E)} \\
&=&\lim_{n}\sum\limits_{k=1}^{n}\left\langle \xi ,T^{-1}\circ \left( T\circ
h_{k}\right) \right\rangle _{M(E)}\left\langle T^{-1}\circ \left( T\circ
h_{k}\right) ,\xi \right\rangle _{M(E)} \\
&=&\lim_{n}\sum\limits_{k=1}^{n}\left\langle \xi ,h_{k}\right\rangle
_{M(E)}\left\langle h_{k},\xi \right\rangle _{M(E)}
\end{eqnarray*}%
for all $\xi \in E$, we deduce that $\{h_{n}\}_{n}$ is a standard frame of
multipliers in $E$. 
\endproof%

\begin{corollary}
If $b(E)$ admits a standard frame of multipliers, then $E$ admits a standard
frame of multipliers.
\end{corollary}

\proof
Let $\{h_{n}\}_{n}$ be a standard frame of multipliers in $b(E)$. By Theorem
3.12 there is an invertible element $T\in L(b(E))$ such that $\{T\circ
h_{n}\}_{n}$ is a standard normalised frame of multipliers in $b(E)$. Since $%
\widetilde{T\circ h_{n}}=\widetilde{T}\circ \widetilde{h_{n}}$ ($\widetilde{%
T\circ h_{n}}$ denotes the extension of $T\circ h_{n}$ to an element in $%
M(E) $ ) for each positive integer $n$, by Corollary 3.9, $\{\widetilde{T}%
\circ \widetilde{h_{n}}\}_{n}$ is a standard normalised frame of multipliers
in $E, $ and then, by Theorem 3.12, $\{\widetilde{h_{n}}\}_{n}$ is a
standard frame of multipliers in $E$, since $\widetilde{T}$ is an invertible
element in $b(L(E))$.%
\endproof%

\section{Dual frames}

\begin{lemma}
Let $E$ be a countably generated Hilbert $A$-module in $M(E)$, let $%
\{h_{n}\}_{n}$ be a standard frame of multipliers in $E$, and let $S$ be an
invertible positive element in $b(L(E))$. Then $\{S\circ h_{n}\}_{n}$ is a
standard frame of multipliers in $E$.
\end{lemma}

\proof%
Let $\xi \in E$. Since 
\begin{equation*}
\sum\limits_{k=1}^{n}\left\langle \xi ,S\circ h_{k}\right\rangle
_{M(E)}\left\langle S\circ h_{k},\xi \right\rangle
_{M(E)}=\sum\limits_{k=1}^{n}\left\langle S^{\ast }\left( \xi \right)
,h_{k}\right\rangle _{M(E)}\left\langle h_{k},S^{\ast }\left( \xi \right)
\right\rangle _{M(E)}
\end{equation*}%
and since $\{h_{n}\}_{n}$ is a standard frame of multipliers in $E$, $%
\sum\limits_{n}\left\langle \xi ,S\circ h_{n}\right\rangle
_{M(E)}\left\langle S\circ h_{n},\xi \right\rangle _{M(E)}$ converges in $A$%
. Moreover, 
\begin{equation*}
C\left\langle S^{\ast }(\xi ),S^{\ast }(\xi )\right\rangle _{E}\leq
\sum\limits_{n}\left\langle \xi ,S\circ h_{n}\right\rangle
_{M(E)}\left\langle S\circ h_{n},\xi \right\rangle _{M(E)}\leq D\left\langle
S^{\ast }(\xi ),S^{\ast }(\xi )\right\rangle _{E}.
\end{equation*}%
Since $S$ is an invertible positive element in $b(L(E))$, $S^{\ast }$ is an
invertible element in $b(L(E))$ and then 
\begin{equation*}
\left\Vert S^{-1}\right\Vert _{\infty }^{-2}\left\langle \xi ,\xi
\right\rangle _{E}\leq \left\langle S^{\ast }(\xi ),S^{\ast }(\xi
)\right\rangle _{E}\leq \left\Vert S\right\Vert _{\infty }^{2}\left\langle
\xi ,\xi \right\rangle _{E}.
\end{equation*}%
From these facts we conclude that $\{S\circ h_{n}\}_{n}$ is a standard frame
of multipliers in $E$. 
\endproof%

\begin{remark}
Let $E$ be a countably generated Hilbert $A$-module in $M(E)$, let $%
\{h_{n}\}_{n}$ be a standard frame of multipliers in $E.$ If $\theta $ is
the frame transform of $\{h_{n}\}_{n},$ then, by Theorem 3.11, $\left(
\theta ^{\ast }\circ \theta \right) ^{-1}$ is an invertible positive element
in $b(L(E))$ and by Lemma 4.1, $\{\left( \theta ^{\ast }\circ \theta \right)
^{-1}\circ h_{n}\}_{n}$ is a standard frame of multipliers in $E$. Let $\xi
\in E$. Since $\{\left( \theta ^{\ast }\circ \theta \right) ^{-\frac{1}{2}%
}\circ h_{n}\}_{n}$ is a standard normalised frame of multipliers in $E$\ (
the proof of Theorem 3.12) we have 
\begin{equation*}
\left( \theta ^{\ast }\circ \theta \right) ^{-\frac{1}{2}}\left( \xi \right)
=\sum\limits_{n}\left( \theta ^{\ast }\circ \theta \right) ^{-\frac{1}{2}%
}\circ h_{n}\cdot \left\langle \left( \theta ^{\ast }\circ \theta \right) ^{-%
\frac{1}{2}}\circ h_{n},\left( \theta ^{\ast }\circ \theta \right) ^{-\frac{1%
}{2}}\left( \xi \right) \right\rangle _{M(E)}.
\end{equation*}%
From this fact and Theorem 3.7, we obtain 
\begin{eqnarray*}
\xi &=&\left( \theta ^{\ast }\circ \theta \right) ^{\frac{1}{2}}\left(
\left( \theta ^{\ast }\circ \theta \right) ^{-\frac{1}{2}}\left( \xi \right)
\right) \\
&=&\left( \theta ^{\ast }\circ \theta \right) ^{\frac{1}{2}}\left(
\sum\limits_{n}\left( \theta ^{\ast }\circ \theta \right) ^{-\frac{1}{2}%
}\circ h_{n}\cdot \left\langle \left( \theta ^{\ast }\circ \theta \right) ^{-%
\frac{1}{2}}\circ h_{n},\left( \theta ^{\ast }\circ \theta \right) ^{-\frac{1%
}{2}}\left( \xi \right) \right\rangle _{M(E)}\right) \\
&=&\left( \theta ^{\ast }\circ \theta \right) ^{\frac{1}{2}}\left(
\sum\limits_{n}\left( \theta ^{\ast }\circ \theta \right) ^{-\frac{1}{2}%
}\circ h_{n}\cdot \left\langle \left( \theta ^{\ast }\circ \theta \right)
^{-1}\circ h_{n},\xi \right\rangle _{M(E)}\right) \\
&=&\left( \theta ^{\ast }\circ \theta \right) ^{\frac{1}{2}}\left(
\lim\limits_{n}\sum\limits_{k=1}^{n}\left( \theta ^{\ast }\circ \theta
\right) ^{-\frac{1}{2}}\circ h_{k}\cdot \left\langle \left( \theta ^{\ast
}\circ \theta \right) ^{-1}\circ h_{k},\xi \right\rangle _{M(E)}\right) \\
&=&\lim\limits_{n}\sum\limits_{k=1}^{n}h_{k}\cdot \left\langle \left( \theta
^{\ast }\circ \theta \right) ^{-1}\circ h_{k},\xi \right\rangle _{M(E)}.
\end{eqnarray*}%
This shows that $\sum\limits_{n}h_{n}\cdot \left\langle \left( \theta ^{\ast
}\circ \theta \right) ^{-1}\circ h_{n},\xi \right\rangle _{M(E)}$ converges
to $\xi .$
\end{remark}

\begin{definition}
Let $E$ be a countably generated Hilbert $A$-module in $M(E)$, let $%
\{h_{n}\}_{n}$ be a standard frame of multipliers in $E$. We say that a
standard frame of multipliers $\{t_{n}\}_{n}$ in $E$ is a dual frame of
multipliers of $\{h_{n}\}_{n}$ if $\sum\limits_{n}h_{n}\cdot \left\langle
t_{n},\xi \right\rangle _{M(E)}$ converges in $E$ for all $\xi \in E$ and
moreover, 
\begin{equation*}
\xi =\sum\limits_{n}h_{n}\cdot \left\langle t_{n},\xi \right\rangle _{M(E)}.
\end{equation*}%
The standard frame of multipliers $\{\left( \theta ^{\ast }\circ \theta
\right) ^{-1}\circ h_{n}\}_{n}$ in $E$, where $\theta $ is the frame
transform of $\{h_{n}\}_{n}$, is called the canonical dual frame of
multipliers of $\{h_{n}\}_{n}$ and $\left( \theta ^{\ast }\circ \theta
\right) ^{-1}$ is said the frame operator of the standard frame of
multipliers $\{h_{n}\}_{n}$.
\end{definition}

\begin{theorem}
( Reconstruction formula) Let $E$ be a countably generated Hilbert $A$%
-module in $M(E)$ and let $\{h_{n}\}_{n}$ be a standard frame of multipliers
in $E$. Then there is a unique invertible positive element $S\in b(L\left(
E\right) )$ such that $\sum\limits_{n}h_{n}\cdot \left\langle S\circ
h_{n},\xi \right\rangle _{M(E)}$ converges in $E$ and moreover, 
\begin{equation*}
\xi =\sum\limits_{n}h_{n}\cdot \left\langle S\circ h_{n},\xi \right\rangle
_{M(E)}
\end{equation*}%
for all $\xi \in E.$
\end{theorem}

\proof
By Theorem 3.12 there is an invertible element $T\in b(L(E))$ such that $%
\{T\circ h_{n}\}_{n}$ is a standard normalized frame of multipliers in $E$.
By Theorem 3.7, for each $\xi \in E,$ $\sum\limits_{n}T\circ h_{n}\cdot
\left\langle T\circ h_{n},\xi \right\rangle _{M(E)}$ converges in $E$ for
each $\xi \in E$, and moreover, 
\begin{equation*}
\xi =\sum\limits_{n}T\circ h_{n}\cdot \left\langle T\circ h_{n},\xi
\right\rangle _{M(E)}.
\end{equation*}%
Then 
\begin{eqnarray*}
&&\overline{p}_{E}\left( \xi -\sum\limits_{k=1}^{n}h_{k}\cdot \left\langle
S\circ h_{k},\xi \right\rangle _{M(E)}\right) \\
&=&\overline{p}_{E}\left( \xi -\sum\limits_{k=1}^{n}h_{k}\cdot \left\langle
T\circ h_{k},T(\xi )\right\rangle _{M(E)}\right) \\
&=&\overline{p}_{E}\left( T^{-1}\left( T(\xi )-\sum\limits_{k=1}^{n}T\circ
h_{k}\cdot \left\langle T\circ h_{k},T(\xi )\right\rangle _{M(E)}\right)
\right) \\
&\leq &\widetilde{p}_{L(E)}(T^{-1})\overline{p}_{E}\left( T(\xi
)-\sum\limits_{k=1}^{n}T\circ h_{k}\cdot \left\langle T\circ h_{k},T(\xi
)\right\rangle _{M(E)}\right) \\
&\leq &\left\Vert T^{-1}\right\Vert _{\infty }\overline{p}_{E}\left( T(\xi
)-\sum\limits_{k=1}^{n}T\circ h_{k}\cdot \left\langle T\circ h_{k},T(\xi
)\right\rangle _{M(E)}\right)
\end{eqnarray*}%
for all $\xi \in E$, for all $p\in S(A)$ and for all positive integer $n$.
From this fact and taking into account that 
\begin{equation*}
T(\xi )=\sum\limits_{n}T\circ h_{n}\cdot \left\langle T\circ h_{n},T(\xi
)\right\rangle _{M(E)}
\end{equation*}%
for all $\xi \in E$, we deduce that $\sum\limits_{n}h_{n}\cdot \left\langle
\left( T^{\ast }\circ T\right) \circ h_{n},\xi \right\rangle _{M(E)}$
converges in $E$, and moreover, 
\begin{equation*}
\xi =\sum\limits_{n}h_{n}\cdot \left\langle \left( T^{\ast }\circ T\right)
\circ h_{n},\xi \right\rangle _{M(E)}.
\end{equation*}%
Let $S=T^{\ast }\circ T.$ Clearly, $S$ is a positive invertible element in $%
b(L(E)).$ Moreover, $\sum\limits_{n}h_{n}\cdot \left\langle S\circ h_{n},\xi
\right\rangle _{M(E)}$ converges in $E$, and $\xi =\sum\limits_{n}h_{n}\cdot
\left\langle S\circ h_{n},\xi \right\rangle _{M(E)}$ for all $\xi \in E.$

To show that $S$ is unique with the above properties, suppose that there is
two positive invertible elements $S_{1}$ and $S_{2}$ in $b(L(E))$ such that
for each $\xi \in E$, $\sum\limits_{n}h_{n}\cdot \left\langle S_{1}\circ
h_{n},\xi \right\rangle _{M(E)}$ and $\sum\limits_{n}h_{n}\cdot \left\langle
S_{2}\circ h_{n},\xi \right\rangle _{M(E)}$ converge in $E$, and 
\begin{equation*}
\xi =\sum\limits_{n}h_{n}\cdot \left\langle S_{1}\circ h_{n},\xi
\right\rangle _{M(E)}=\sum\limits_{n}h_{n}\cdot \left\langle S_{2}\circ
h_{n},\xi \right\rangle _{M(E)}.
\end{equation*}%
Then 
\begin{eqnarray*}
\xi &=&\sum\limits_{n}h_{n}\cdot \left\langle S_{1}\circ h_{n},\xi
\right\rangle _{M(E)}=\sum\limits_{n}h_{n}\cdot \left\langle S_{1}\circ
S_{2}^{-1}\circ S_{2}\circ h_{n},\xi \right\rangle _{M(E)} \\
&=&\sum\limits_{n}h_{n}\cdot \left\langle S_{2}\circ h_{n},\left(
S_{2}^{-1}\circ S_{1}\right) \left( \xi \right) \right\rangle _{M(E)}=\left(
S_{2}^{-1}\circ S_{1}\right) \left( \xi \right)
\end{eqnarray*}%
for all $\xi \in E$. This implies that $S_{1}=S_{2}$ and the uniqueness is
proved. 
\endproof%

\begin{remark}
The dual frame of multipliers of a given standard frame of multipliers is
unique.
\end{remark}

\begin{proposition}
Let $E$ be a countably generated Hilbert $A$-module in $M(E)$ and let $%
\{h_{n}\}_{n}$ and $\{t_{n}\}_{n}$ be two standard frames of multipliers in $%
E$ with the frame transforms $\theta _{1}$ and $\theta _{2}$. Then these
frames of multipliers are duals to each other if and only if $\theta
_{1}^{\ast }\circ \theta _{2}=$id$_{E}.$
\end{proposition}

\proof%
First we suppose that the standard frames of multipliers $\{h_{n}\}_{n}$ and 
$\{t_{n}\}_{n}$ are duals to each other. Then 
\begin{eqnarray*}
\left\langle \left( \theta _{1}^{\ast }\circ \theta _{2}\right) \left( \xi
\right) ,\eta \right\rangle _{E} &=&\left\langle \theta _{2}\left( \xi
\right) ,\theta _{1}\left( \eta \right) \right\rangle _{H_{A}} \\
&=&\dsum\limits_{n}\left\langle \xi ,t_{n}\right\rangle _{M(E)}\left\langle
h_{n},\eta \right\rangle _{M(E)}=\lim_{n}\dsum\limits_{k=1}^{n}\left\langle
\xi ,t_{k}\right\rangle _{M(E)}\left\langle h_{k},\eta \right\rangle _{M(E)}
\\
&=&\lim_{n}\left\langle \dsum\limits_{k=1}^{n}h_{k}\cdot \left\langle
t_{k},\xi \right\rangle _{M(E)},\eta \right\rangle _{M(E)} \\
&=&\left\langle \lim_{n}\left( \dsum\limits_{k=1}^{n}h_{k}\cdot \left\langle
t_{k},\xi \right\rangle _{M(E)}\right) ,\eta \right\rangle _{M(E)} \\
&=&\left\langle \xi ,\eta \right\rangle _{E}
\end{eqnarray*}%
for all $\xi ,\eta \in E,$ and so $\theta _{1}^{\ast }\circ \theta _{2}=$id$%
_{E}$.

Conversely, suppose that $\theta _{1}^{\ast }\circ \theta _{2}=$id$_{E}$.
Let $\xi .$ From

$\overline{p}_{E}\left( \dsum\limits_{k=n+1}^{\infty }h_{k}\cdot
\left\langle t_{k},\xi \right\rangle _{M(E)}\right) $

$=\sup \{p\left( \left\langle \dsum\limits_{k=n+1}^{\infty }h_{k}\cdot
\left\langle t_{k},\xi \right\rangle _{M(E)},\eta \right\rangle _{E}\right)
;\eta \in E,\overline{p}_{E}\left( \eta \right) \leq 1\}$

$=\sup \{p\left( \dsum\limits_{k=n+1}^{\infty }\left\langle \xi
,t_{k}\right\rangle _{M(E)}\left\langle h_{k},\eta \right\rangle
_{M(E)}\right) ;\eta \in E\overline{p}_{E}\left( \eta \right) \leq 1\}$

$\leq \sup \{p\left( \dsum\limits_{k=n+1}^{\infty }\left\langle \xi
,t_{k}\right\rangle _{M(E)}\left\langle t_{k},\xi \right\rangle
_{M(E)}\right) ^{\frac{1}{2}}p\left( \dsum\limits_{k=n+1}^{\infty
}\left\langle \eta ,h_{k}\right\rangle _{M(E)}\left\langle h_{k},\eta
\right\rangle _{M(E)}\right) ^{\frac{1}{2}};\ \ $

$\eta \in E,\overline{p}_{E}\left( \eta \right) \leq 1\}$

$\leq p\left( \dsum\limits_{k=n+1}^{\infty }\left\langle \xi
,t_{k}\right\rangle _{M(E)}\left\langle t_{k},\xi \right\rangle
_{M(E)}\right) ^{\frac{1}{2}}\sup \{p\left( \dsum\limits_{n}\left\langle
\eta ,h_{n}\right\rangle _{M(E)}\left\langle h_{n},\eta \right\rangle
_{M(E)}\right) ^{\frac{1}{2}};$

$\ \eta \in E,\overline{p}_{E}\left( \eta \right) \leq 1\}$

$\leq p\left( \dsum\limits_{k=n+1}^{\infty }\left\langle \xi
,t_{k}\right\rangle _{M(E)}\left\langle t_{k},\xi \right\rangle
_{M(E)}\right) ^{\frac{1}{2}}\sup \{p\left( D_{1}\left\langle \eta ,\eta
\right\rangle _{E}\right) ^{\frac{1}{2}};\eta \in E,\overline{p}_{E}\left(
\eta \right) \leq 1\}$

$\leq D_{1}p\left( \dsum\limits_{k=n+1}^{\infty }\left\langle \xi
,t_{k}\right\rangle _{M(E)}\left\langle t_{k},\xi \right\rangle
_{M(E)}\right) ^{\frac{1}{2}}$

$\ $for all $p\in S(A)$ and taking into account that $\dsum\limits_{n}\left%
\langle \xi ,t_{n}\right\rangle _{M(E)}\left\langle t_{n},\xi \right\rangle
_{M(E)}$ converges in $A$, we deduce that $\dsum\limits_{n}h_{n}\cdot
\left\langle t_{n},\xi \right\rangle _{M(E)}$ converges in $A$.

Since

\begin{equation*}
\left\langle \xi ,\eta \right\rangle _{E}=\left\langle \left( \theta
_{1}^{\ast }\circ \theta _{2}\right) \left( \xi \right) ,\eta \right\rangle
_{E}=\dsum\limits_{n}\left\langle \xi ,t_{n}\right\rangle
_{M(E)}\left\langle h_{n},\eta \right\rangle _{M(E)}
\end{equation*}%
for all $\eta \in E$, we have 
\begin{eqnarray*}
&&p\left( \left\langle \xi -\dsum\limits_{k=1}^{n}h_{k}\cdot \left\langle
t_{k},\xi \right\rangle _{M(E)},\eta \right\rangle _{E}\right) \\
&=&p\left( \left\langle \xi ,\eta \right\rangle -\left\langle
\dsum\limits_{k=1}^{n}h_{k}\cdot \left\langle t_{k},\xi \right\rangle
_{M(E)},\eta \right\rangle _{E}\right) \\
&=&p\left( \dsum\limits_{n}\left\langle \xi ,t_{n}\right\rangle
_{M(E)}\left\langle h_{n},\eta \right\rangle _{M(E)}-\left\langle
\dsum\limits_{k=1}^{n}h_{k}\cdot \left\langle t_{k},\xi \right\rangle
_{M(E)},\eta \right\rangle _{E}\right) \\
&=&p\left( \dsum\limits_{n}\left\langle \xi ,t_{n}\right\rangle
_{M(E)}\left\langle h_{n},\eta \right\rangle
_{M(E)}-\dsum\limits_{k=1}^{n}\left\langle \xi ,t_{k}\right\rangle
_{M(E)}\left\langle h_{k},\eta \right\rangle _{M(E)}\right) \\
&=&p\left( \dsum\limits_{k=n+1}^{\infty }\left\langle \xi
,t_{k}\right\rangle _{M(E)}\left\langle h_{k},\eta \right\rangle
_{M(E)}\right)
\end{eqnarray*}%
for all positive integer $n$, for all $\eta \in E$ and for all $p\in S(A)$
and since $\dsum\limits_{n}\left\langle \xi ,t_{n}\right\rangle _{M(E)}$ $%
\left\langle h_{n},\eta \right\rangle _{M(E)}$ converges in $A$, 
\begin{equation*}
\lim\limits_{n}\left\langle \xi -\dsum\limits_{k=1}^{n}h_{k}\cdot
\left\langle t_{k},\xi \right\rangle _{M(E)},\eta \right\rangle _{E}=0
\end{equation*}%
for all $\eta \in E.$ From these facts, we deduce that $\dsum%
\limits_{n}h_{n}\cdot \left\langle t_{n},\xi \right\rangle _{M(E)}$
converges to $\xi $ and so the standard frames of multipliers $\{h_{n}\}_{n}$
and $\{t_{n}\}_{n}$ are duals to each other.%
\endproof%

\begin{corollary}
Let $E$ be a countably generated Hilbert $A$-module in $M(E)$. The canonical
bi-dual frame of a standard frame of multipliers $\{h_{n}\}_{n}$ in $E$ is
the frame itself.
\end{corollary}

\proof%
Let $\xi \in E$. If $\theta ^{\prime }$ is the frame transform of the dual
frame of multipliers $\{\left( \theta ^{\ast }\circ \theta \right)
^{-1}\circ h_{n}\}_{n},$ then

\begin{eqnarray*}
\theta ^{\prime }\left( \xi \right) &=&\left( \left\langle \left( \theta
^{\ast }\circ \theta \right) ^{-1}\circ h_{n},\xi \right\rangle \right)
_{n}=\left( \left\langle h_{n},\left( \theta ^{\ast }\circ \theta \right)
^{-1}\left( \xi \right) \right\rangle \right) _{n} \\
&=&\theta \left( \left( \theta ^{\ast }\circ \theta \right) ^{-1}\left( \xi
\right) \right) =\left( \theta \circ \left( \theta ^{\ast }\circ \theta
\right) ^{-1}\right) \left( \xi \right) .
\end{eqnarray*}%
This shows that $\theta ^{\prime }=\theta \circ \left( \theta ^{\ast }\circ
\theta \right) ^{-1}$. A simple calculus show that 
\begin{equation*}
\left( \left( \theta ^{\prime }\right) ^{\ast }\circ \theta ^{\prime
}\right) ^{-1}=\theta ^{\ast }\circ \theta
\end{equation*}%
and then 
\begin{equation*}
\left( \left( \theta ^{\prime }\right) ^{\ast }\circ \theta ^{\prime
}\right) ^{-1}\circ \left( \theta ^{\ast }\circ \theta \right) ^{-1}\circ
h_{n}=h_{n}
\end{equation*}%
for all positive integer $n$. Moreover, 
\begin{equation*}
\left( \theta ^{\prime }\right) ^{\ast }\circ \theta =\text{id}_{E}\text{.}
\end{equation*}%
From these facts and Proposition 4.6, we conclude that the canonical bi-dual
frame of multipliers of the standard frame of multipliers $\{h_{n}\}_{n}$ in 
$E$ is the frame itself.%
\endproof%


\begin{thebibliography}{10}
\bibitem[1]{1} M. Fragoulopoulou, \textit{Topological algebras with
involution, } North-Holland Mathematics Studies, 200. Elsevier Science B.V.,
Amsterdam, 2005.

\bibitem[2]{2} M. Frank, D. R. Larson, \textit{A module frame concept for
Hilbert }$C^{\ast }$\textit{-modules}, in Functional and Harmonic Analysis
of Wavelets (San Antonio, TX, Jan. 1999), A.M.S., Providence, R.I., Contemp.
Math. 247 (2000), 207-233.

\bibitem[3]{3} M. Frank and D. R. Larson,\textit{\ Frames in Hilbert }$%
C^{\ast }$\textit{-modules and }$C^{\ast }$\textit{-algebras}. J. Operator
Theory \textbf{48} (2002), no. 2, 273--314.

\bibitem[4]{4} A. Inoue, \textit{Locally }$C^{\ast }$\textit{-algebras,}
Mem. Faculty Sci. Kyushu Univ. Ser. A,\textit{\ }\textbf{25}(1971), 197-235.%
\textit{\ }

\bibitem[5]{5} M. Joi\c{t}a, \textit{Projections on Hilbert modules over
locally }$C^{\ast }$\textit{-algebras,} Math. Reports, \textbf{4(54)}(2002),
pp.373-378 (2003).

\bibitem[6]{6} M. Joi\c{t}a, \textit{On Hilbert-modules over locally }$%
C^{\ast }$\textit{-algebras. II.} Period. Math. Hungar. \textbf{51} (2005),
no. 1, 27--36.

\bibitem[7]{7} M. Joi\c{t}a, \textit{Hilbert modules over locally }$C^{\ast
} $\textit{-algebras}, University of Bucharest Press, (2006), 150 pg. ISBN
973737128-3.

\bibitem[8]{8} M. Joita, \textit{On multiplier modules of Hilbert modules
over locally }$C^{\ast }$\textit{-algebras} ( submitted), arXiv:0707.1139v1.

\bibitem[9]{9} A. Khosravi and M.\ S. Asgari, \textit{\ Frames and bases in
Hilbert modules over locally }$C^{\ast }$\textit{-algebras}. Int. J. Pure
Appl. Math. \textbf{14} (2004), no. 2, 171--190.

\bibitem[10]{10} G. G. Kasparov, \textit{Hilbert }$C^{\ast }$\textit{%
-modules: Theorem of Stinespring and Voiculescu,} J. Operator Theory \textbf{%
4}(1980), 133-150.

\bibitem[11]{11} E. C. Lance, \textit{Hilbert }$C^{\ast }$\textit{-modules.
A toolkit for operator algebraists, }London Mathematical Society Lecture
Note Series 210, Cambridge University Press, Cambridge 1995.

\bibitem[12]{12} N. C. Phillips, \textit{Inverse limits of }$C^{\ast }$%
\textit{-algebras,} J. Operator Theory, \textbf{19 }(1988), 159-195.

\bibitem[13]{13} I. Raeburn and S. J. Thompson, \textit{Countably generated
Hilbert modules, the Kasparov stabilisation theorem, and frames with Hilbert
modules,} Proc. Amer. Math. Soc. \textbf{131} (2003), no. 5, 1557--1564
\end{thebibliography}
\end{document}